\documentclass[preprint,11pt]{elsarticle}
\usepackage{amssymb}
\usepackage{amsmath}
\usepackage{amsthm}
\usepackage{amsfonts}
\usepackage{bbm}
\usepackage{enumerate}
\usepackage{epic}
\usepackage{setspace}
\numberwithin{equation}{section}
\numberwithin{figure}{section}
\numberwithin{table}{section}

\newtheorem{theorem}{Theorem}[section]

\newtheorem{proposition}{Proposition}[section]
\newproof{pf}{Proof}
\newdefinition{example}{Example}[section]
\newdefinition{remark}{Remark}[section]

\begin{document}

\begin{frontmatter}
\title{Poly-Bernoulli numbers and lonesum matrices}

\author{Hyun Kwang Kim\fnref{author1}}
\ead{hkkim@postech.ac.kr}
\address{Department of Mathematics, Pohang University of Science and
Technology, Pohang 790-784, Korea}

\author{Denis S. Krotov\fnref{author2}}
\ead{krotov@math.nsc.ru}
\address{Sobolev Institute of Mathematics, pr. Akademika Koptyuga 4,
Novosibirsk 630090, Russia}
\address{Mechanics and Mathematics Department, Novosibirsk State
University, Russia}

\author{Joon Yop Lee\corref{corauthor}\fnref{author3}}
\ead{flutelee@postech.ac.kr}
\address{ASARC, KAIST, Daejeon 305-340, Korea}

\cortext[corauthor]{Corresponding author}

\fntext[author1]{Tel \#: +82-54-279-2049, 
 Fax \#: +82-54-279-2799.\\
The work of Hyun Kwang Kim was supported by Basic Research
Program through the National Research Foundation of Korea (NRF)
funded by the Ministry of Education, Science and Technology (grant
\# 2009-0089826).}

\fntext[author2]{Tel \#: +7-383-3634666, Fax \#: +7-383-3332598.\\
The work of Denis S. Krotov was supported by the Federal Target Grant
``Scientific and Educational Personnel of Innovation Russia''
for 2009-2013 (contract No. 02.740.11.0429) and the
Russian Foundation for Basic Research (grant \# 10-01-00424)}

\fntext[author3]{Tel \#: +82-42-350-8111, Fax \#: +82-42-350-8114.\\ 
The work of Joon Yop Lee was supported by Basic Science Research Program through the National Research Foundation of Korea (NRF) funded by the Ministry of Education, Science and Technology (grant \# 2010-0001655).}

\begin{abstract}
A lonesum matrix is a matrix that can be uniquely reconstructed from
its row and column sums. Kaneko defined the poly-Bernoulli numbers $B_m^{(n)}$ by a generating function, and Brewbaker computed the number of binary lonesum $m\times
n$-matrices and showed that this number coincides with the poly-Bernoulli number $B_m^{(-n)}$. We compute the number of $q$-ary lonesum $m\times n$-matrices, and then provide generalized Kaneko's formulas by using the generating function for the number of
$q$-ary lonesum $m\times n$-matrices. In addition, we define two types of $q$-ary lonesum matrices that are composed of strong and weak lonesum matrices, and suggest further researches on lonesum matrices. 
\end{abstract}


\end{frontmatter}

\section{Introduction\label{section1}}

Kaneko \cite{kaneko,kaneko1} defined and studied the poly-Bernoulli numbers $B_m^{(n)}$
of index $n\in\mathbb{Z}$ by the generating function
\begin{equation}\label{equation1}
\frac{Li_n(1-e^{-x})}{1-e^{-x}}=\sum_{m=0}^{\infty}B_m^{(n)}\frac{x^m}{m!}
\end{equation}
where $Li_n(z)$ denotes the formal power series
$\sum_{l=1}^{\infty}\frac{z^l}{l^n}$ (the $n$th polylogarithm when
$n>0$ and the rational function $\big(z\frac
{d}{dz}\big)^{-n}\big(\frac{z}{1-z}\big)$ when $n\leq 0$). By
analyzing this generating function, Kaneko et al. proved formulas for the
poly-Bernoulli numbers of negative indices \cite{kaneko,arakawa}:
\begin{align}
&B_m^{(-n)}=\sum_{l=0}^m(-1)^{l+m}l!S(m,l)(l+1)^n,\label{equation2}\\
&B_m^{(-n)}=\sum_{l=0}^{min(m,n)}(l!)^2S(m+1,l+1)S(n+1,l+1),\label{equation3}\\
&\sum_{n=0}^{\infty}\sum_{m=0}^{\infty}B_m^{(-n)}\frac{x^m}{m!}
\frac{y^n}{n!}=\frac{e^{x+y}}{e^x+e^y-e^{x+y}}.\label{equation4}
\end{align}
Here $S(m,l)$ denotes the Stirling number of the second kind, which
is the number of ways to partition an $m$-element set into $l$
nonempty subsets. S\'anchez-Peregrino also proved the equation (\ref{equation3}) by a much simpler way \cite{sanchez,sanchez1}.
The motivation of this research is combinatorial
interpretations of Kaneko's formulas.

A \emph{binary matrix} is a matrix each of whose entries is either 0
or 1 and a \emph{lonesum matrix} is a matrix that can be uniquely
reconstructed from its row and column sums. Ryser \cite{ryser}
proved that a binary matrix is a lonesum matrix if and only if each
of its $2\times 2$-submatrices is not
$$\begin{pmatrix}1& 0\\0& 1\end{pmatrix}\ \text{and}\ \begin{pmatrix}0& 1\\1& 0\end{pmatrix}.$$
Brewbaker exploited this result and
computed the number of binary lonesum $m\times n$-matrices by the
principle of inclusion and exclusion \cite{brewbaker1,brewbaker}. This number
coincides with the right hand side of the formula (\ref{equation2}).

We investigate properties of binary lonesum matrices in Section
\ref{section2}. First of all, we present the second proof of Ryser's
theorem. Brewbaker may have known this proof although he did not
explicitly state it in \cite{brewbaker1, brewbaker}. The main idea of the second
proof is to partition the set of binary lonesum
$m\times n$-matrices. Computing the number of lonesum matrices in each set and
summing them up, Brewbaker obtained the right hand side of the
formula (\ref{equation3}). Therefore we may regard that Kaneko's
formulas (\ref{equation2}) and (\ref{equation3}) are computations of
the number of binary lonesum $m\times n$-matrices, on one hand, by
using the principle of inclusion and exclusion and, on the other
hand, by partitioning the set of binary lonesum $m\times n$-matrices.

A \emph{$q$-ary matrix} is a matrix each of whose entries is in $\{0,1,\ldots,q-1\}$. We consider $q$-ary lonesum matrices in
Section \ref{section3}. There are two types of lonesum matrices in
$q$-ary matrices, namely, strong and weak lonesum matrices (see the
definitions in Section \ref{section3}). Throughout Subsection
\ref{subsection3.1} (resp. \ref{subsection3.2}), a lonesum matrix
means a strong (resp. weak) lonesum matrix.

We devote Subsection \ref{subsection3.1} to strong lonesum matrices.
We first generalize the second proof of Ryser's theorem to $q$-ary
matrices. It turns out that every minimal $q$-ary non-lonesum matrix
is of size $2\times 2$. We next compute the number of ternary lonesum 
$m\times n$-matrices and generalize this technique to $q$-ary
lonesum matrices. We also give the formula for the number of
$q$-ary symmetric lonesum $n\times n$-matrices. We finally compute
the generating function for the number of $q$-ary lonesum matrices
and generalize Kaneko's formulas (\ref{equation1}) and
(\ref{equation4}).

We dedicate Subsection \ref{subsection3.2} to weak lonesum matrices.
After studying properties of weak lonesum matrix, we show that 
if $q\ge 5$ then the number of $q$-ary weak non-lonesum matrices is infinite, and we
construct some $q$-ary weak non-lonesum matrices when $q\in\{3,4\}$.
We also suggest an open problem related with weak lonesum matrices.

\section{Binary lonesum matrices\label{section2}}

A binary matrix is a matrix each of whose entries is either $0$ or
$1$. Throughout this section, every matrix is a binary matrix unless
we specify otherwise. A lonesum matrix is a matrix that can be
uniquely reconstructed it from its row and column sums. For example,
a $3\times 3$ matrix with rows $(1,1,0)$, $(1,0,0)$, and $(1,1,0)$
is a lonesum matrix because of the unique reconstruction:
$$
\begin{matrix}
\cline{1-3}
\multicolumn{1}{|c}*&*&*&\multicolumn{1}{|c}2\\
\multicolumn{1}{|c}*&*&*&\multicolumn{1}{|c}1\\
\multicolumn{1}{|c}*&*&*&\multicolumn{1}{|c}3\\
\cline{1-3}
3&2&1&
\end{matrix}
\longrightarrow
\begin{pmatrix}
1&1&0\\
1&0&0\\
1&1&1
\end{pmatrix}$$

Brewbaker proved most results of this section \cite{brewbaker1, brewbaker}. For
the convenience of the reader and the sake of completeness, we give
a few properties of lonesum matrices with short proofs, which are
useful in our discussion of $q$-ary lonesum matrices in Section
\ref{section3}.

The criterion to distinguish lonesum matrices from non-lonesum ones
is a theorem proved by Ryser \cite{ryser}.
\begin{theorem}[Ryser]\label{ryser}
A binary matrix is a lonesum matrix if and only if each of its
$2\times 2$ submatrices is not
$$\begin{pmatrix}1 & 0\\0& 1\end{pmatrix}\ \text{and}\
\begin{pmatrix}0 & 1\\1& 0\end{pmatrix}.$$
\end{theorem}
\begin{pf}
See \cite{ryser, brualdi1, haber, ryser1, ryser2}.
\end{pf}
\noindent Two matrices $M$ and $M^{\prime}$ are said to be
\emph{equivalent} if $M$ is changed into $M^{\prime}$ by a row
or column permutations, that is, there are permutation
matrices $P$ and $P^{\prime}$ such that $M^{\prime}=PAP^{\prime}$. A
\emph{forbidden matrix} is a non-lonesum matrix each of whose proper
submatrices is a lonesum matrix, that is, a minimal non-lonesum
matrix. Evidently, every matrix that contains a forbidden matrix is
not a lonesum matrix. Hence we focus on forbidden matrices that are
not equivalent to each other. Using these terminologies, we can
rephrase Theorem \ref{ryser} in the form that
$\begin{pmatrix}1& 0\\0& 1\end{pmatrix}$ is the unique forbidden
matrix.

Let $M$ be a lonesum $m \times n$-matrix. Identifying a binary
vector with the set of its nonzero coordinates, we regard a row of
$M$ as a subset of $\{1,2, \ldots,n\}$. Then Theorem 2.1 states that
$M$ is a lonesum matrix if and only if the rows of $M$ form a chain
in the inclusion lattice formed by the subsets of
$\{1,2,\ldots,n\}$. Using this fact, Brewbaker proved that the right
hand side of Kaneko's formula (\ref{equation2}) gives the number of
lonesum $m \times n$-matrices \cite{brewbaker1, brewbaker}.

Using that the rows of a lonesum matrix form a chain, we consider
another proof of Theorem \ref{ryser}. This proof provides the main
tool for constructing lonesum matrices and computing the number of
lonesum $m\times n$-matrices.

{\bf The second proof of Theorem \ref{ryser}.}  A
\emph{stair} matrix is a matrix whose $i$th row is
$(1^{r_i},0^{n-r_i})$ with $r_1 \ge r_2 \ge \cdots\geq r_m$.
Figure \ref{bstandard} is an example of stair matrix.
\begin{figure}[h]
$$
\begin{matrix}
&4&3&10&2&5&8&1&9&6&7\\
\cline{2-11}
2&\multicolumn{1}{|c}1&1&1&1&1&1&1&1&0&\multicolumn{1}{c|}0\\
1&\multicolumn{1}{|c}1&1&1&1&1&1&0&0&0&\multicolumn{1}{c|}0\\
5&\multicolumn{1}{|c}1&1&1&1&1&1&0&0&0&\multicolumn{1}{c|}0\\
6&\multicolumn{1}{|c}1&1&1&1&1&1&0&0&0&\multicolumn{1}{c|}0\\
7&\multicolumn{1}{|c}1&1&1&0&0&0&0&0&0&\multicolumn{1}{c|}0\\
11&\multicolumn{1}{|c}1&1&1&0&0&0&0&0&0&\multicolumn{1}{c|}0\\
3&\multicolumn{1}{|c}1&0&0&0&0&0&0&0&0&\multicolumn{1}{c|}0\\
4&\multicolumn{1}{|c}1&0&0&0&0&0&0&0&0&\multicolumn{1}{c|}0\\
4&\multicolumn{1}{|c}1&0&0&0&0&0&0&0&0&\multicolumn{1}{c|}0\\
9&\multicolumn{1}{|c}1&0&0&0&0&0&0&0&0&\multicolumn{1}{c|}0\\
\cline{2-11}
\end{matrix}$$
\caption{A stair matrix\label{bstandard}}
\end{figure}
It is sufficient to
prove that a matrix has no submatrix equivalent to
$\begin{pmatrix}1&0\\0&1\end{pmatrix}$ if and only if it is
equivalent to a stair matrix.

The necessity is clear because every stair matrix is a lonesum
matrix. For the sufficiency, let $M$ be an $m \times n$-matrix
that has no submatrix equivalent to
$\begin{pmatrix}1&0\\0&1\end{pmatrix}$. Permuting the rows of $M$
suitably, we assume that $r_1 \ge r_2 \ge \cdots \ge r_m$ where
$r_i$ denotes the number of $1$s in the $i$th row of $M$. By an
appropriate column permutation, we further assume that the first row
of $M$ is $(1^{r_1},0^{n-r_1})$. Then the last $n-r_1$ entries of
the second row of $M$ should be all zeros. Now a suitable column
permutation changes the second row of $M$ into $(1^{r_2},0^{n-r_2})$
with $r_1 \ge r_2$. If we continue this process, then $M$ becomes
equivalent to a stair matrix. \mbox{}\hfill\qed

Let $M$ be a lonesum $m\times n$-matrix. The second proof of
Theorem \ref{ryser} implies that $M$ is equivalent to a stair
matrix, called the standard form of $M$. Hence an ordered partition
pair
\begin{equation}
\big((A_0,A_1,\ldots,A_j),(B_0,B_1,\ldots,B_j)\big)
\label{partitionpair}\end{equation} where $$\biguplus_{i=0}^j
A_i=\{1,2,\ldots,n\},\, \biguplus_{i=0}^jB_i=\{1,2,\ldots,k\}$$ with
$|A_0|\ge 0$, $|B_0|\ge 0$ and $|A_i|\ge 1$ $|B_i|\ge 1$ for $i\ge 1$
determines the positions of $1$ in $M$. Here $\biguplus$ denotes the
disjoint union.  For example, the matrix in Figure \ref{bstandard} satisfies
$$\begin{cases}
(A_0,A_1,A_2,A_3,A_4)=
\big(\emptyset,\{2\},\{1,5,6\},\{7,11\},
\{3,4,9,10\}\big)\\
(B_0,B_1,B_2,B_3,B_4)=
\big(\{4\},\{3,10\},\{2,5,8\},\{1,9\},\{6,7\}\big)
\end{cases}.$$
This means that the $(a,b)$-entry of $M$ is $1$ if
and only if $(a,b)\in A_i\times\big(\biguplus_{h=0}^{j-i}B_h\big)$
for some $i$. We call $j+1$ the number of \emph{stairs} in $M$.
Since $|A_0|\ge 0$ and $|B_0|\ge 0$, the construction of
(\ref{partitionpair}) is equivalent to that of
$$\Big(\big(\{0\}\cup A_0, A_1,\ldots,A_j\big),\big(\{0\}\cup
B_0,B_1,\ldots,B_j\big)\Big).$$ Therefore for a fixed $j$ the number
of ordered partition pairs of the form (\ref{partitionpair}) is
\begin{equation}\label{lonesum1}(j!)^2S(m+1,j+1)S(n+1,j+1).\end{equation}
Summing (\ref{lonesum1}) for $0\leq j\leq min(m,n)$ yields the
number of binary lonesum $m\times n$-matrices
$$\sum_{j=0}^{min(m,n)}(j!)^2S(m+1,j+1)S(n+1,j+1),$$
which coincides with the formula (\ref{equation3}).

This proves that Kaneko's formula (\ref{equation3}) is nothing but
computation of the number of lonesum $m \times n$-matrices by
partitioning the set of lonesum $m\times n$-matrices according to
the standard forms of its elements. Brewbaker may have known this
fact \cite{brewbaker1,brewbaker}. We include this proof since we can directly
generalize the idea of this proof to $q$-ary lonesum matrices in
Section \ref{section3}.

\begin{remark} Using quantum algebra, Launois \cite{launois,launois1} proved that
the number of lonesum $m \times n$-matrices equals the number of
permutations $\sigma$ on $\{1,\,2,\ldots,\,m+n \}$ satisfying $-n
\le \sigma(i)-i \le m $ for $1\leq i\leq m+n$. We can construct a
one-to-one correspondence between such permutations and lonesum $m \times n$-matrices as follows:

\begin{enumerate}[{Step} 1.]
\item Let $M=(M_{ji})$ be a lonesum $m\times n$-matrix with exactly
$k$ different nonzero rows and $k$ different nonzero columns. The reason of denoting the entries of $M$ by $M_{ji}$ is that we need to determine a position of a {\em rook} (defined later) in $2$-dimensional space by an $M_{ji}$. Assigning an order to the columns sums of $M$, we define a $k$-tuple
$(C_1,C_2,\ldots ,C_k)$ of subsets of $\{1,2,\ldots ,n\}$ by the rule: $i\in
C_a$ if and only if the $i$th column has the $a$th column sum. Similarly, we define
a $k$-tuple $(R_1,R_2,\ldots ,R_k)$ of subsets of $\{1,2,\ldots ,m\}$, which corresponds to the row sums of $M$. For
example,
$$
\begin{matrix}
&C_2&C_3&&C_3&C_1&C_2&C_2&C_3\\
\cline{2-9}
R_2&\multicolumn{1}{|c}1&1&0&1&0&1&1&\multicolumn{1}{c|}1\\
R_1&\multicolumn{1}{|c}0&1&0&1&0&0&0&\multicolumn{1}{c|}1\\
R_3&\multicolumn{1}{|c}1&1&0&1&1&1&1&\multicolumn{1}{c|}1\\
R_2&\multicolumn{1}{|c}1&1&0&1&0&1&1&\multicolumn{1}{c|}1\\
\cline{2-9}
\end{matrix}
$$
Using the standard form of lonesum matrix,
it is easy to
understand that the mapping
$$ M \to \big((C_1,C_2,\ldots ,C_k),(R_1,R_2,\ldots ,R_k)\big)$$ defines a one-to-one correspondence 
between the lonesum $m\times n$-matrices that have exactly $k$ different column sums
(equivalently, exactly $k$ different row sums)
and the pairs of $k$-tuples composed of mutually disjoint
nonempty subsets of $\{1,2,\ldots ,n\}$ and $\{1,2,\ldots ,m\}$, respectively.

For convenience, we define two partition functions
$$\begin{cases}C:\{1,2,\ldots ,n\}\to \{0,1,\ldots ,k\}\\
R:\{1,2,\ldots ,m\}\to \{0,1,\ldots ,k\}\end{cases}$$
 by
$$\begin{cases}
C(i)=a&\text{for}\ i\in C_a\\
R(j)=b&\text{for}\ j\in R_b\end{cases}$$
where
$$\begin{cases}
C_0 = \{1,2,\ldots ,n\} \setminus \underset{a=1}{\overset{k}{\bigcup}}C_a\\
R_0 = \{1,2,\ldots ,m\} \setminus \underset{b=1}{\overset{k}{\bigcup}}R_b
\end{cases}.$$
The inverse mapping is defined by
$$ M_{ji} = \begin{cases}
1 & \text{for}\ $C(i)+R(j)>k$\\
0& \text{ otherwise}.
\end{cases}.$$

\item Considering Step 1, we establish a one-to-one correspondence between the
 lonesum $m\times n$-matrices and the permutations
$\sigma$ on $\{1,2,\ldots ,m+n\}$ that satisfy 
\begin{equation}\label{rcondition}
-n \le \sigma(i)-i \le m
\end{equation}
for $i\in \{1,2,\ldots ,m+n\}$.

We start with a permutation matrix $P$ that corresponds to $\sigma$.
We call the ones of the matrix $P$ \emph{rooks}.
Assume that $P$ satisfies the condition (\ref{rcondition}).
Let $k$ be the number of rooks in the top-left $n\times m$-submatrix of $P$, and $$(i_1,j_1),(i_2,j_2),\ldots, (i_k,j_k)$$
be the positions of these rooks
where $j_1<j_2<\ldots<j_k$.

Next, we define a partition function
$$C:\{1,2,\ldots ,n\}\to\{0,1,\ldots ,k\}.$$ For $a\in\{1,2,\ldots,k\}$ we let $C(i_a)=a$
$$\begin{array}{c@{}|@{}c@{}c@{}c@{}c@{}|@{}c@{}c@{}c@{}c@{}c@{}c@{}c@{}|}
\hline
 2-&-&-&1& & &0&0&0&0&0&0\\
 3-&-&-&-&1& & &0&0&0&0&0\\
   & & & & & & & &0&0&0&0\\
   & & & & & & & & &0&0&0\\
 1-&1& & & & & & & & &0&0\\
 \hphantom{-}&\hphantom{-}&\hphantom{-}&\hphantom{-}&\hphantom{-}&\hphantom{-}&\hphantom{-}&\hphantom{-}&\hphantom{-}&\hphantom{-}&\hphantom{-}&
\end{array}$$
and if there is a rook in the positoin $(i,m+1)$ for some $i\in\{1,2,\ldots,n\}$, then we define $C(i)=0$
$$\begin{array}{c@{}|@{}c@{}c@{}c@{}c@{}|@{}c@{}c@{}c@{}c@{}c@{}c@{}c@{}|}
\hline
   & & & & & &0&0&0&0&0&0\\
   & & & & & & &0&0&0&0&0\\
 0-&-&-&-&-&1& & &0&0&0&0\\
 \hphantom{-}&\hphantom{-}&\hphantom{-}&\hphantom{-}&\hphantom{-}&\hphantom{-}&\hphantom{-}&\hphantom{-}&\hphantom{-}&\hphantom{-}&\hphantom{-}&
\end{array}$$
For any other rooks in the positions $(i,j)$ for $i\in\{1,2,\ldots,n\}$, we apply
the rule $C(i)=C(j-m-1)$
$$\begin{array}{c@{}|@{}c@{}c@{}c@{}c@{}|@{}c@{}c@{}c@{}c@{}c@{}c@{}c@{}|}
\hline
   & & & & & &0&0&0&0&0&0\\
   & & & &1&-&-&\backslash&0&0&0&0\\
   & & & & & & &|&0&0&0&0\\
   & & & & & & &1& &0&0&0\\
   & & & & & & & & & &0&0\\
 \hphantom{-}&\hphantom{-}&\hphantom{-}&\hphantom{-}&\hphantom{-}&\hphantom{-}&\hphantom{-}&\hphantom{-}&\hphantom{-}&\hphantom{-}&\hphantom{-}&
\end{array}$$
Since $j-m-1<i$, the function $C$ is well-defined. 
Remark that the sets $C_1,C_2,\ldots,C_k$ are nonempty while the set
$C_0$ can be empty. 

Similarly, the $m$ rooks with the $m$ smallest values $j$ in their position $(i,j)$ of the matrix $P$ define the function
$$R:\{1,2\ldots ,m\}\to\{0,1,\ldots ,k\}.$$ 
For this, we rotate the matrix $P$
to the angle of $180^\circ$ and then apply the same algorithm, by exchanging the roles of $m$ and $n$.

\item Now we describe the inverse transformation
$$\big((C_1,C_2,\ldots,C_k),(R_1,R_2,\ldots,R_k)\big)\to M.$$
Assume that we have two partitions $(C_0,C_1,\ldots ,C_k)$
and $(R_0,R_1,\ldots ,R_k)$ of $\{1,2,\ldots,n\}$ and $\{1,2,\ldots,m\}$, respectively, where only $C_0$ and $R_0$ can be empty.
For $i\in\{0,1,\ldots,k\}$ let $$C_i=\{a_{i,1},a_{i,2},\ldots,a_{i,|C_i|}\}$$ satisfying
$a_{i,1}<a_{i,2}\ldots<a_{i,|C_i|}$.
The positions of $|C_i|$ rooks that correspond to the set $C_i$ are
$$(a_{i,1},c_i),(a_{i,2},a_{i,1}+m+1),\ldots,(a_{i,|C_i|},\,a_{i,|C_i|-1}+m+1)$$
where $c_0=m+1$, $1 \le c_1<c_2 \ldots  <c_k \le m$,
and exact values of $c_1,c_2,\ldots,c_k$ are still unknown.
So, we already know the positions of the 
$n-k$ rooks in the top-right $n\times n$-submatrix of $P$,
and the first coordinates and the order (from the left to the right)
of the positions of the $k$ rooks in the top-left $n \times m$-submatrix of $P$.

Similarly (in a symmetrical way),
the partition $(R_0,R_1,\ldots ,R_k)$ defines the positions of the 
$n-k$ rooks in the bottom-left $m\times m$ submatrix of $P$
and the first coordinates and the order
of the positions of the $k$ rooks in the lower-right $m \times n$ submatrix of $P$.
Since no two rooks can not be in the same column of $P$,
this information uniquely identifies the permutation matrix $P$.
\end{enumerate}

Lov$\acute{a}$sz gave a
slightly weaker version of this correspondence \cite{lovasz}.

\end{remark}

\section{$q$-ary lonesum matrices\label{section3}}

A $q$-ary matrix is a matrix each of whose entries is in
$\{0,1,\ldots,q-1\}$. There are two types of lonesum matrices for
$q$-ary matrices. For a $q$-ary vector $\mathbf{v}=(v_1,\ldots,v_n)$
the \emph{structure vector} of $\mathbf{v}$ is
$(a_0,a_1,\ldots,a_{q-1})$ where $a_j$ is the number of entries in
$\mathbf{v}$ such that $v_i=j$.  A \emph{strong} (resp. \emph{weak})
lonesum matrix is a matrix that can be uniquely reconstructed
from its row and column sums (resp. structure vectors). For binary
matrices, the definition of strong lonesum matrix and that of weak
one are identical, however, this is not true for nonbinary matrices.
\begin{example} A $3 \times 3$-matrix with the rows
$(0,1,0)$, $(1,2,1)$, and $(0,1,0)$ is a strong lonesum matrix
because we can uniquely reconstruct it from its row and column sums:
$$
\begin{matrix}
\cline{1-3}
\multicolumn{1}{|c}*&*&\multicolumn{1}{c|}*&1\\
\multicolumn{1}{|c}*&*&\multicolumn{1}{c|}*&4\\
\multicolumn{1}{|c}*&*&\multicolumn{1}{c|}*&1\\
\cline{1-3}
1&4&1&
\end{matrix}
\longrightarrow\ 
\begin{pmatrix}
0&1&0\\
1&2&1\\
0&1&0
\end{pmatrix}
.$$

A $3 \times 3$-matrix with the rows $(0,1,0)$, $(1,2,1)$, and
$(0,1,1)$ is a weak lonesum matrix because we can uniquely
reconstruct it from its row and column structure vectors:
\begin{equation*}
\begin{matrix}
\cline{1-3}
\multicolumn{1}{|c}*&*&\multicolumn{1}{c|}*&(2,1,0)\\
\multicolumn{1}{|c}*&*&\multicolumn{1}{c|}*&(0,2,1)\\
\multicolumn{1}{|c}*&*&\multicolumn{1}{c|}*&(1,2,0)\\
\cline{1-3}
(2,1,0)&(0,2,1)&(1,2,0)&\end{matrix}\  
\longrightarrow\ \begin{pmatrix}0&1&0\\1&2&1\\0&1&1\end{pmatrix}.\end{equation*}
Notice that this matrix is not a strong lonesum matrix. In fact, we
can construct two different matrices with the row sums $1$, $4$, $2$
and the column sums $1$, $4$, $2$:
$$
\begin{matrix}
\cline{1-3}
\multicolumn{1}{|c}0&1&\multicolumn{1}{c|}0&1\\
\multicolumn{1}{|c}1&2&\multicolumn{1}{c|}1&4\\
\multicolumn{1}{|c}0&1&\multicolumn{1}{c|}1&2\\
\cline{1-3}
1&4&2&
\end{matrix}
\ \text{and}\hspace{3mm}
\begin{matrix}
\cline{1-3}
\multicolumn{1}{|c}0&1&\multicolumn{1}{c|}0&1\\
\multicolumn{1}{|c}1&1&\multicolumn{1}{c|}2&4\\
\multicolumn{1}{|c}0&2&\multicolumn{1}{c|}0&2\\
\cline{1-3}
1&4&2&
\end{matrix}$$
\end{example}

\subsection{Strong lonesum matrices\label{subsection3.1}}
It follows from Theorem \ref{ryser} that every binary forbidden
matrix is equivalent to $\begin{pmatrix}1&0\\0&1\end{pmatrix}$. This
means that every $2\times 2$-submatrix of a binary lonesum matrix is
equivalent to one of
\begin{gather*}
\begin{pmatrix}1&1\\1&1\end{pmatrix}
,\,\begin{pmatrix}1&b\\c&0\end{pmatrix},\,\text{and}\ 
\begin{pmatrix}0&0\\0&0\end{pmatrix}
\end{gather*}
where $\{b,c\}\subseteq\{0,1\}$. By using this criterion, we have computed
the number of binary lonesum $m\times n$-matrices. We consider a
similar criterion for $q$-ary lonesum matrices and compute the
number of $q$-ary lonesum $m\times n$-matrices by exploiting this
new criterion.

\subsubsection{The criterion for $q$-ary lonesum matrices}
We first consider $q$-ary lonesum $2\times 2$-matrices.
Let $$M=\begin{pmatrix}a&b\\
c&d\end{pmatrix}\ \text{and}\ M(\alpha)=\begin{pmatrix}a-\alpha&b+\alpha\\
c+\alpha&d-\alpha\end{pmatrix}$$ for $\alpha\in\mathbb{Z}$ be $q$-ary
matrices. If $M$ is a lonesum matrix, then $M$ and $M(\alpha)$ have
the same row and column sums if and only if $\alpha=0$. Without loss
of generality, we assume that $max\{a,b,c,d\}=a$.

Suppose that $a=q-1$. If either $b=q-1$ or $c=q-1$, then $\alpha=0$.
We assume that $b,\,c<q-1$. Then $\alpha=0$ if and only if $d=0$.
Hence $\alpha=0$ if and only if $M$ is one of
\begin{equation*}\label{qcriterion1}
\begin{pmatrix}q-1&q-1\\c&d\end{pmatrix},\,
\begin{pmatrix}q-1&b\\q-1&d\end{pmatrix},\,
\begin{pmatrix}q-1&b\\c&0\end{pmatrix}.
\end{equation*}
Suppose that $a<q-1$. Similarly, $\alpha=0$ if and only if $M$ is
one of
\begin{equation*}\label{qcriterion2}
\begin{pmatrix}a&b\\0&0\end{pmatrix}\ \text{and}\
\begin{pmatrix}a&0\\c&0\end{pmatrix}.
\end{equation*}
Therefore $\alpha=0$ if and only if $M$ is one of
\begin{equation}\label{forbidden2}
\begin{pmatrix}q-1&q-1\\c&d\end{pmatrix},\,
\begin{pmatrix}q-1&b\\q-1&d\end{pmatrix},\,
\begin{pmatrix}q-1&b\\c&0\end{pmatrix},\,
\begin{pmatrix}a&b\\0&0\end{pmatrix},\,
\begin{pmatrix}a&0\\c&0\end{pmatrix}.
\end{equation}

This implies that if a $q$-ary matrix is a lonesum matrix, then each
of its $2\times 2$-submatrices is equivalent to one of the matrices
(\ref{forbidden2}). By using this, we show that this is also a
sufficient condition for a $q$-ary matrix to be a lonesum matrix.
Since we can easily generalize the criterion for ternary lonesum
matrices to that for $q$-ary ones, we first consider the criterion
for ternary lonesum matrices. Recall that a \emph{ternary} matrix is
a matrix each of whose entries is in $\{0,1,2\}$.
\begin{theorem}\label{discriminantt}
A ternary matrix is a lonesum matrix if and only if each of its
$2\times 2$-submatrices is equivalent to one of
\begin{equation}\label{forbidden}
\begin{pmatrix}2&2\\c&d\end{pmatrix},\,
\begin{pmatrix}2&b\\2&d\end{pmatrix},\,
\begin{pmatrix}2&b\\c&0\end{pmatrix},\,
\begin{pmatrix}a&b\\0&0\end{pmatrix},\,
\begin{pmatrix}a&0\\c&0\end{pmatrix}
\end{equation}
where $a,\,b,\,c,\,d\in\{0,1,2\}$.
\end{theorem}
\noindent Without proof, Brualdi stated this theorem \cite{brualdi} .
\begin{pf}
Similar to the second proof of Theorem \ref{ryser}, we need only
show that if every $2\times 2$-submatrix of a ternary matrix is
equivalent to one of the matrices (\ref{forbidden}), then it is a
lonesum matrix.

To prove this, we use induction on the size of matrix. Let $M$ be an
ternary $m\times n$-matrix each of whose $2\times 2$-submatrices is
equivalent to one of the matrices (\ref{forbidden}) and $M_{ij}$
be the $(i,j)$-entry of $M$. We denote by
$(r_{i,0},r_{i,1},r_{i,2})$ (resp. $(c_{j,0},c_{j,1},c_{j,2})$) the
structure vector of the $i$th row (resp. $j$th column). Permuting
the rows and columns of $M$, we assume that the row structure
vectors of $M$ satisfy
$$\begin{cases}
r_{i,2}\ge r_{i+1,2}\\
r_{i,2}+r_{i,1}\ge r_{i+1,2}+r_{i+1,2}\\
r_{i,2}+r_{i,1}+r_{i,0}=r_{i+1,2}+r_{i+1,1}+r_{i+1,0}
\end{cases}$$
and the column structure vectors of $M$ satisfy
$$\begin{cases}
c_{j,2}\ge c_{j+1,2}\\
c_{j,2}+c_{j,1}\ge c_{j+1,2}+c_{j+1,1}\\
c_{j,2}+c_{j,1}+c_{j,0}=c_{j+1,2}+c_{j+1,1}+c_{j+1,0}
\end{cases}.$$
Every $2\times 2$-submatrix of $M$ is equivalent to one of the
matrices (\ref{forbidden}), thus, by a simple reasoning,
$$M_{i_1j_1}\ge M_{i_2j_2}\ \text{if and only if}\ i_1\leq i_2\ \text{and}\ j_1\leq
j_2.$$ We call this the \emph{standard form} of $M$. Figure \ref{tstandard} is an example of standard form.
\begin{figure}[h]
$$\begin{pmatrix}
2&2&2&2&2&2&2&1&1&0&0\\
2&2&2&2&2&2&2&0&0&0&0\\
2&2&2&1&1&0&0&0&0&0&0\\
2&2&2&0&0&0&0&0&0&0&0\\
2&2&2&0&0&0&0&0&0&0&0\\
2&2&2&0&0&0&0&0&0&0&0\\
2&1&0&0&0&0&0&0&0&0&0\\
2&1&0&0&0&0&0&0&0&0&0\\
2&1&0&0&0&0&0&0&0&0&0\\
2&0&0&0&0&0&0&0&0&0&0
\end{pmatrix}$$
\caption{A ternary standard form\label{tstandard}}
\end{figure}

Suppose that either $m=1$ or $n=1$. Then $M$ is evidently a lonesum
matrix.

Suppose that $m\ge 2$ and $n\ge 2$. We denote by $M_i$ the $i$th row
of $M$. Defining $\Big\lceil\frac{|M_{l+1}|}{2}\Big\rceil=-1$, we
let $l$ be the number such that
$$\bigg\lceil\frac{|M_1|}{2}\bigg\rceil=\bigg\lceil\frac{|M_2|}{2}\bigg\rceil=
\cdots=\bigg\lceil\frac{|M_l|}{2}\bigg\rceil>\bigg\lceil\frac{|M_{l+1}|}{2}\bigg\rceil$$
where $|M_i|$ is the sum of entries in $M_i$.

Suppose that $l=1$. Each $2\times 2$-submatrix of $M$ is equivalent
to one of the matrices (\ref{forbidden}), thus
$$M=\begin{pmatrix}
2&\cdots&2&M_{1,k+1}&\cdots&M_{1,n}\\
M_{2,1}&\cdots&M_{2,k}&0&\cdots&0\\
\cdots&\cdots&\cdots&\cdots&\cdots&\cdots\\
M_{n,1}&\cdots&M_{n,k}&0&\cdots&0
\end{pmatrix}$$ for some $k$. Since we know the column sums of $M$, we can
uniquely reconstruct the first row of $M$. Eliminating the first row
of $M$ yields an $(m-1)\times n$-matrix $M'$ each of whose $2\times
2$-submatrices is one of the matrices (\ref{forbidden}). If we
use induction on the size of $M'$, then we can reconstruct $M'$ from
its row and column sums. Since we can reconstruct both $M'$ and the
first row of $M$ from the row and column sums of $M$, we can
reconstruct $M$ from its row and column sums.

Suppose that $l\ge 2$. Similar to the case $l=1$, the $l\times n$
matrix $M^{\prime}$ formed by the first $l$ rows of $M$ is
$$M^{\prime}=\begin{pmatrix}
2&\cdots&2&M_{1,k^{\prime}+1}&0&\cdots&0\\
2&\cdots&2&M_{2,k^{\prime}+1}&0&\cdots&0\\
\cdots&\cdots&\cdots&\cdots&\cdots&\cdots&\cdots\\
2&\cdots&2&M_{l,k^{\prime}+1}&0&\cdots&0
\end{pmatrix}$$ for some $k^{\prime}$.
Since we know the row sums of $M$, we can uniquely reconstruct
$M^{\prime}$ from the row sums of $M$. Eliminating the first $l$
rows of $M$ changes $M$ into an $(m-l)\times n$-matrix each of whose
$2\times 2$-submatrices is one of the matrices (\ref{forbidden}).
Similar to the case $l=1$, we can reconstruct $M$ from its row and
column sums.
\end{pf}
\begin{remark}\mbox{}\hfill
\label{remark3.2}\begin{enumerate}[1.]
\item Theorem \ref{discriminantt} implies that every ternary forbidden
matrix is equivalent to one of
\begin{gather*}\begin{pmatrix}2&b\\c&2\end{pmatrix},\,
\begin{pmatrix}2&b\\c&1\end{pmatrix},\,\begin{pmatrix}1&b\\c&1\end{pmatrix}
\end{gather*}
where $\{b,c\}\subseteq\{0,1\}$.

\item Analyzing the proof of Theorem \ref{discriminantt}, we can
construct a ternary lonesum matrix $M$ as follows.
\begin{enumerate}[{Step} 1.]
\item Permuting the rows and columns of $M$, we assume that
$M_{ij}\ge M_{i^{\prime}j^{\prime}}$ if and only if $i\ge
i^{\prime}$ and $j\ge j^{\prime}$.

\item Form stairs in $M$ with $2$.

\item Determine positions of both $0$s and $1$s by considering the
$2\times 2$-matrices (\ref{forbidden}) allowed in $M$.
\end{enumerate}
\end{enumerate}
\end{remark}

In fact, we can apply the proof of Theorem \ref{discriminantt} to
$q$-ary matrices by substituting the role of $q-1$ (resp.
$\{0,1,\ldots,q-2\}$) for that of $2$ (resp. $\{0,1\}$). This
application yields the criterion for $q$-ary lonesum matrices.
\begin{theorem}\label{discriminantq}
A $q$-ary matrix is a lonesum matrix if and only if each of its
$2\times 2$-submatrices is equivalent to one of
\begin{equation}\label{forbiddenq}
\begin{pmatrix}q-1&q-1\\c&d\end{pmatrix},\,
\begin{pmatrix}q-1&b\\q-1&d\end{pmatrix},\,
\begin{pmatrix}q-1&b\\c&0\end{pmatrix},\,
\begin{pmatrix}a&b\\0&0\end{pmatrix},\,
\begin{pmatrix}a&0\\c&0\end{pmatrix}
\end{equation}
where $\{a,b,c,d\}\subseteq\{0,1,\ldots,q-2\}$.
\end{theorem}
\begin{remark}
Theorem \ref{discriminantq} says that every non-lonesum matrix
contains a $2\times 2$-matrix not equivalent to one of the matrices
(\ref{forbidden2}). Hence each $q$-ary forbidden matrix is of
size $2\times 2$.
\end{remark}

\subsubsection{The number of $q$-ary lonesum
$m\times n$-matrices}Using the second remark of Remark \ref{remark3.2}, we can
construct ternary lonesum matrices, and comparing Theorems
\ref{discriminantt} and \ref{discriminantq}, we can apply
construction of ternary lonesum matrices to that of $q$-ary ones.
Hence, for simplicity, we first consider the case of ternary lonesum
matrices.

Let $M$ be a ternary lonesum $m\times n$-matrix. According to the
criterion for ternary lonesum matrices, the positions $2$ in $M$ is determined by a partition pair
$$\Big(\big(\{0\}\cup A_0,A_1,\ldots,A_j\big),\big(\{0\}\cup B_0,B_1,\ldots,B_j\big)\Big)$$
where
$$\begin{cases}
\underset{i=0}{\overset{j}{\biguplus}}A_i=\{1,2,\ldots,m\}\\
\underset{i=0}{\overset{j}{\biguplus}}B_i=\{1,2,\ldots,n\}\end{cases}.$$
We define the \emph{$i$th
block} of $M$ to be
$$\big\{(a,b)\ \big|\ (a,b)\in A_i\times B_{j+1-i}\big\}$$
where $1\leq i\leq j$. For example, the blocks of the matrix in Figure \ref{tstair} are
$$\begin{cases}
A_1\times B_3=\{2,6\}\times\{1,6,7,10\}\\
A_2\times B_2=\{1,5,7,11\}\times\{2,5,8,9\}\\
A_3\times B_1=\{3,4,9,10\}\times\{3,11\}
\end{cases}.$$
\begin{figure}[h]
$$
\begin{matrix}
&4&3&11&2&5&8&9&1&6&7&10\\
\cline{2-12}
2&\multicolumn{1}{|c}2&2&2&2&2&2&2&\multicolumn{1}{|c}{\mathbf{0}}&\mathbf{0}&\mathbf{0}&\multicolumn{1}{c|}{\mathbf{0}}\\
6&\multicolumn{1}{|c}2&2&2&2&2&2&2&\multicolumn{1}{|c}{\mathbf{0}}&\mathbf{0}&\mathbf{0}&\multicolumn{1}{c|}{\mathbf{0}}\\
\cline{5-12}
1&\multicolumn{1}{|c}2&2&2&\multicolumn{1}{|c}{\mathbf{0}}&\mathbf{0}&\mathbf{0}&\multicolumn{1}{c|}{\mathbf{0}}&0&0&0&\multicolumn{1}{c|}0\\
5&\multicolumn{1}{|c}2&2&2&\multicolumn{1}{|c}{\mathbf{0}}&\mathbf{0}&\mathbf{0}&\multicolumn{1}{c|}{\mathbf{0}}&0&0&0&\multicolumn{1}{c|}0\\
7&\multicolumn{1}{|c}2&2&2&\multicolumn{1}{|c}{\mathbf{0}}&\mathbf{0}&\mathbf{0}&\multicolumn{1}{c|}{\mathbf{0}}&0&0&0&\multicolumn{1}{c|}0\\
11&\multicolumn{1}{|c}2&2&2&\multicolumn{1}{|c}{\mathbf{0}}&\mathbf{0}&\mathbf{0}&\multicolumn{1}{c|}{\mathbf{0}}&0&0&0&\multicolumn{1}{c|}0\\
\cline{3-8}
3&\multicolumn{1}{|c}2&\multicolumn{1}{|c}{\mathbf{0}}&\multicolumn{1}{c|}{\mathbf{0}}&0&0&0&0&0&0&0&\multicolumn{1}{c|}0\\
4&\multicolumn{1}{|c}2&\multicolumn{1}{|c}{\mathbf{0}}&\multicolumn{1}{c|}{\mathbf{0}}&0&0&0&0&0&0&0&\multicolumn{1}{c|}0\\
9&\multicolumn{1}{|c}2&\multicolumn{1}{|c}{\mathbf{0}}&\multicolumn{1}{c|}{\mathbf{0}}&0&0&0&0&0&0&0&\multicolumn{1}{c|}0\\
10&\multicolumn{1}{|c}2&\multicolumn{1}{|c}{\mathbf{0}}&\multicolumn{1}{c|}{\mathbf{0}}&0&0&0&0&0&0&0&\multicolumn{1}{c|}0\\
\cline{2-12}
\end{matrix}
$$
\caption{Blocks of a ternary matrix\label{tstair}}
\end{figure}The criterion for ternary lonesum
matrices implies that $1$s can be in the blocks of $M$ only.

For an $r\times s$-block the position set of $1$s should be, by the criterion for ternary lonesum matrices,
$$\big\{(x_1,y),\ldots,(x_t,y)\big\}\ \text{or}\
\big\{(x,y_1),\ldots,(x,y_u)\big\}.$$ There are four cases:
$$\begin{cases}
\text{Case 1: No $1$ in the block}\\
\text{Case 2: $t=1$ or $u=1$}\\
\text{Case 3: $t\ge 2$}\\
\text{Case 4: $u\ge 2$}\end{cases}.$$ By a simple computation, the number of position sets of
$1$s in each case is the following: $$\begin{cases}
\text{Case 1: $1$}\\
\text{Case 2: $r s$}\\
\text{Case 3: $s\overset{r}{\underset{l=2}{\sum}}{r\choose l}=s(2^{r}-r-1)$}\\
\text{Case 4:
$r\overset{s}{\underset{l=2}{\sum}}{s\choose
l}=r(2^{s}-s-1)$}\end{cases}.$$ Summing the
numbers in Cases 1--4 yields the number of position sets of $1$s in the $r\times s$-block
$$f_3(r,s)=1+r s+r(2^{s}-s-1)+s(2^r-r-1).$$

For a ternary $m\times n$-matrix the number of ways to form blocks of sizes
$m_1\times n_j,m_1\times n_{j-1},\ldots,m_j\times n_1$
is $${m\choose
m_0,m_1,\ldots,m_j}{n\choose
n_0,n_1,\ldots,n_j}.$$
Therefore, denoting 
$$\mathcal{S}_l^j=\Bigg\{ (l_0,l_1,\ldots,l_j)\in\Bbb{Z}^{j+1}\
\Bigg|\ \sum_{i=0}^jl_i=l,\,l_0\ge 0,\,l_i\ge 1\ \big(i\in\{1,2,\ldots,j\}\big)\Bigg\}$$ for $j\ge 1$ and $l\ge 1$, if we consider all of the standard forms and possible positions of $1$s in blocks, then we obtain the number of
ternary lonesum $m\times n$-matrices.
\begin{theorem}\label{numbert}
The number of ternary lonesum $m\times n$-matrices
$B_m^{(-n)}(3)$ is
\begin{equation}\label{ternarylonesum}
1+\sum_{j=1}^{min(m,n)}\sum_{{(m_0,m_1,\ldots,m_j)\in\mathcal{S}_m^j}\atop
{(n_0,n_1,\ldots,n_j)\in\mathcal{S}_n^j}} {m\choose
m_0,m_1,\ldots,m_j}{n\choose
n_0,n_1,\ldots,n_j}\prod_{i=1}^jf_3(m_{i},n_{j+1-i}).\end{equation}
\end{theorem}

We now consider the case of $q$-ary lonesum matrices. Using Theorem
\ref{discriminantq} and generalizing the technique for ternary
lonesum matrices, we can compute the number of $q$-ary
lonesum $m\times n$-matrices. For this computation, we need only change
$f_3(r,s)$ in the formula (\ref{ternarylonesum}) into
$$f_q(r,s)=1+(q-2)rs+r\big((q-1)^s-(q-2)s-1\big)+s\big((q-1)^r-(q-2)r-1\big),$$
which is the number of position sets of $0$s, $1$s, \ldots, $(q-2)$s in an
$r \times s$-block.
\begin{theorem}\label{numberq}
The number of $q$-ary lonesum $m\times n$-matrices $B_n^{(-k)}(q)$
is
\begin{equation}\label{qarylonesum}
1+\sum_{j=1}^{min(m,n)}\sum_{{(m_0,m_1,\ldots,m_j)\in\mathcal{S}_m^j}\atop
{(n_0,n_1,\ldots,n_j)\in\mathcal{S}_n^j}} {m\choose
m_0,m_1,\ldots,m_j}{n\choose
n_0,n_1,\ldots,n_j}\prod_{i=1}^jf_q(m_{i},n_{j+1-i}).\end{equation}
\end{theorem}

\subsubsection{Symmetric lonesum matrices} We can uniquely
reconstruct a symmetric lonesum matrix from its row or column sums.
Hence both construction of symmetric lonesum matrices and
computation of the number of those are simpler than those of
ordinary ones.

To construct a $q$-ary symmetric lonesum $n\times n$-matrix, we
need only choose an ordered partition $\big(\{0\}\cup A_0,A_1,\ldots,A_j\big)$ of $\{0,1,\ldots,n\}$ instead of an ordered partition pair. In addition, we may pair
$\big(\{0\}\cup A_0,A_1\big)$ and $\big(A_{2i},A_{2i+1}\big)$ for $i\ge
1$ to form blocks. By the criterion for $q$-ary lonesum matrices,
if the parity of $j$ is even, then the block $\big(A_j,A_j\big)$
forms a diagonal matrix with at most one nonzero entry. Therefore,
if we determine positions of $0$s, $1$s, \ldots, $(q-2)$s by the criterion
for $q$-ary lonesum matrices, then we obtain the number of $q$-ary
symmetric lonesum $n\times n$-matrices.
\begin{theorem}
The number of $q$-ary symmetric lonesum $n\times n$-matrices
$B_n(q)$ is
\begin{align*}
&1+\sum_{j=1}^{n}\sum_{(n_0,n_1,\ldots,n_j)\in\mathcal{S}_n^j}
{n\choose n_0,n_1,\ldots,n_j}
\bigg(\prod_{i=1}^{\lfloor\frac{j}{2}\rfloor}
f_q(n_{2i-1},n_{2i})\bigg)\\
&\hspace{9em}
\bigg(1+(q-2)(n-\sum_{i=0}^{2\lfloor\frac{j}{2}\rfloor}n_{i}
)\bigg).
\end{align*}
\end{theorem}
For $n\in\{1,2,3,4,5\}$ the number of $q$-ary symmetric lonesum $n\times n$-matrices is given in Table \ref{table1}.
\begin{table}[h]
\begin{tabular}{|c||c|c|c|}
\hline$n$&$B_n(q)$&$B_n(2)$&$B_n(3)$\\ \hline\hline
$1$       &$q$&$2$&$3$  \\\hline
$2$        &$2q^2+2q-6$&$6$&$18$\\\hline
$3$    &$9q^3-12q^2+12q-22$&$26$& $149$  \\\hline
$4$     &$16q^4+72q^3-312q^2+392q-218$&$150$& $1390$\\\hline
$5$  & $25q^5+160q^4+400q^3-3180q^2+4920q-2598$&$1082$&$13377$ \\\hline
\end{tabular}
\vspace{1em}
\caption{The number of symmetric $q$-ary lonesum matrices}\label{table1}
\end{table}\\
It is known that $B_n(2)$ is the number of necklaces of partitions of $n+2$ labeled beads \cite{sloane}, however, we do not know any combinatorial meaning of the numbers $B_n(q)$ for $q\ge 3$. For example, the terms $3, 18, 149, 1390, 13377$ for $B_n(3)$ do not match any sequence in \cite{sloane}. We think that it is a fascinating task to find a combinatorial object that explains a combinatorial meaning of the numbers $B_n(q)$ for $q\ge 3$.

\subsubsection{Generating functions for the number of lonesum
matrices and generalizations of Kaneko's formulas} Kaneko
\cite{kaneko, kaneko1} defined the poly-Bernoulli numbers $B^{(n)}_m$ by the
generating function (\ref{equation1}) and found out the formula
(\ref{equation4}). By $q$-ary lonesum matrices, we have generalized
the poly-Bernoulli numbers of negative indices. We calculate
generating functions for the number of $q$-ary lonesum matrices.
This provides generalizations of Kaneko's formulas (\ref{equation1})
and (\ref{equation4}).

The exponential generating function for the number of ternary
lonesum matrices is
\begin{align*}
&\sum_{m=0}^{\infty}\sum_{n=0}^{\infty}B_m^{(-n)}(3)\frac{x^m}{m!}\frac{y^n}{n!}\\
=&\sum_{m=0}^{\infty}\sum_{n=0}^{\infty}
\Bigg(1+\sum_{j=1}^{min(m,n)}\sum_{{(m_0,m_1,\ldots,m_j)\in\mathcal{S}_m^j}\atop
{(n_0,n_1,\ldots,n_j)\in\mathcal{S}_n^j}} {m\choose
m_0,m_1,\ldots,m_j}{n\choose
n_0,n_1,\ldots,n_j}\\
&\hspace{3em}\prod_{i=1}^jf_3(m_i,n_{j+1-i})\Bigg)\frac{x^m}{m!}\frac{y^n}{n!}\\
=&\sum_{m=0}^{\infty}\frac{x^m}{m!}\sum_{n=0}^{\infty}\frac{y^n}{n!}\sum_{l=0}^{\infty}
\Bigg(\sum_{r=1}^{\infty}\sum_{s=1}^{\infty}f_3(r,s)\frac{x^r}{r!}\frac{y^s}{s!}\Bigg)^l\\
=&\frac{e^{x+y}}{1-\sum_{r=1}^{\infty}\sum_{s=1}^{\infty}f_3(r,s)\frac{x^r}{r!}\frac{y^s}{s!}}.
\end{align*}
Let $F_3(x,y)=\overset{\infty}{\underset{r=1}{\sum}}
\overset{\infty}{\underset{s=1}{\sum}}f_3(r,s)\frac{x^r}{r!}\frac{y^s}{s!}$.
By
$$
\overset{\infty}{\underset{r=0}{\sum}}\frac{x^r}{r!}=e^x\
\text{and}\
\overset{\infty}{\underset{r=0}{\sum}}r\frac{x^r}{r!}=xe^x,$$ we
gain
$$F_3(x,y)=1-e^x-e^y+(1-x-y-xy+xe^y+ye^x)e^{x+y}.$$ Therefore the
exponential generating function for the number of ternary lonesum
matrices is
$$\sum_{m=0}^{\infty}\sum_{n=0}^{\infty}B_m^{(-n)}(3)\frac{x^m}{m!}\frac{y^n}{n!}
=\frac{e^{x+y}}{1-F_3(x,y)}.$$

Similarly, if we use
\begin{align*}
F_q(x,y)&=\sum_{r=1}^{\infty}\sum_{s=1}^{\infty}f_q(r,s)\frac{x^r}{r!}\frac{y^s}{s!}\\
&=1-e^x-e^y+\big(1-x-y-(q-2)xy+xe^{(q-2)y}+ye^{(q-2)x}\big)e^{x+y},
\end{align*}
then we obtain the exponential generating function for the number of
$q$-ary lonesum matrices
\begin{equation}\label{generalization1}
\sum_{m=0}^{\infty}\sum_{n=0}^{\infty}B_m^{(-n)}(q)\frac{x^m}{m!}\frac{y^n}{n!}=\frac{e^{x+y}}{1-F_q(x,y)}.
\end{equation}
The formula (\ref{generalization1}) is a generalization of Kaneko's
formula (\ref{equation4}).

For the case of symmetric lonesum matrices, if we apply the
computation of the exponential generating function for the number of
$q$-ary lonesum matrices and use
$$\sum_{n=0}^{\infty}\big(1+(q-2)n\big)\frac{x^{n}}{n!}=\big(1+x(q-2)\big)e^x,$$
then we obtain the exponential generating function for the number of
symmetric $q$-ary lonesum matrices
$$\sum_{n=0}^{\infty}B_n(q)\frac{x^n}{n!}=\frac{\big(1+x(q-2)\big)e^{2x}}{1-F_q(x,x)}.$$

 Now we
consider a generalization of (\ref{equation4}). By the definition of
$F_q(x,y)$,
\begin{align*}
F_q(x,y)^l&=\Big\{(1-e^x)+(-1+e^x-xe^x)e^y\\
&\hspace{1em}+\big(-1-(q-2)x+e^{(q-2)x}\big)e^xye^y+xe^xe^{(q-1)y}\Big\}^l\\
          &=\sum_{l_1+l_2+l_3+l_4=l\atop l_1,\,l_2,\,l_3,\,l_4\ge
          0}{l\choose
          l_1,l_2,l_3,l_4}(1-e^x)^{l_1}(-1+e^x-xe^x)^{l_2}x^{l_4}e^{(l_3+l_4)x}
         \\
          &\hspace{7em}\big(-1-(q-2)x+e^{(q-2)x}\big)^{l_3}
         \sum_{m=0}^{\infty}\frac{\big(l_2+l_3+(q-1)l_4\big)^my^{l_3+m}}
          {m!}.
          \end{align*}
Hence the generating function (\ref{generalization1}) becomes
\begin{align}\label{generalization2}
&\sum_{m=0}^{\infty}\sum_{n=0}^{\infty}B_m^{(-n)}(q)\frac{x^m}{m!}\frac{y^n}{n!}\\
 =&\sum_{n=0}^{\infty}\sum_{0\leq l_3\leq n\atop
l_1,l_2,l_4\ge 0 }{ l_1+l_2+l_3+l_4\choose
l_1,l_2,l_3,l_4}(1-e^x)^{l_1}(-1+e^x-xe^x)^{l_2}x^{l_4}e^{(l_3+l_4+1)x}
          &\nonumber\\
&\hspace{7em}\big(-1-(q-2)x+e^{(q-2)x}\big)^{l_3}
\frac{\big(l_2+l_3+(q-1)l_4\big)^{n-l_3}
n!}{(n-l_3)!}\frac{y^n}{n!}. &\nonumber
\end{align}
Computing the coefficient of $\frac{x^m}{m!}$ in
(\ref{generalization2}) generalizes the generating function
(\ref{equation1}):
\begin{align*}&\sum_{n=0}^{\infty}B_n^{(-k)}(q)\frac{x^n}{n!}\\
=&\sum_{0\leq l_3\leq k\atop l_1,\,l_2,\,l_4,\ge 0
}l_3!{l_1+l_2+l_3+l_4\choose l_1,l_2,l_3,l_4}{n\choose
l_3}(1-e^x)^{l_1}(-1+e^x-xe^x)^{l_2}x^{l_4}e^{(l_3+l_4+1)x}
         \\
&\hspace{5.5em} \big(-1-(q-2)x+e^{(q-2)x}\big)^{l_3}
\big(l_2+l_3+(q-1)l_4\big)^{n-l_3}
\end{align*}

\subsection{Weak lonesum matrices\label{subsection3.2}}

We have found out all the forbidden matrices for $q$-ary strong lonesum matrices and these are $2\times 2$-matrices. From this, we naturally wonder if there is a finite number of forbidden matrices for $q$-ary weak lonesum matrices. After studying properties of weak lonesum matrices, we show that 
if $q\ge 5$ then the number of forbidden matrices for $q$-ary weak lonesum matrix is infinite. We also construct some  nontrivial forbidden ternary and quarternary matrices.
Recall that a \emph{quaternary} matrix means a $4$-ary matrix.

Let $M$ be a $q$-ary matrix. A sequence $S=(M_{i_1,j_1},M_{i_2,j_2},\ldots,M_{i_k,j_k})$ 
of mutually different entries of $M$ is called a \emph{path} or an $ab$-\emph{path}
if it satisfies the following conditions:
\begin{enumerate}[{Condition} 1.]
\item For each $l\in\{1,2,\ldots,k-1\}$ either $i_l=i_{l+1}$ or $j_l=j_{l+1}$ is true.

\item For each $l\in\{1,2,\ldots,k-2\}$ both $\big|\{i_l,i_{l+1},i_{l+2}\}\big|\ge 2$ and $\big|\{j_l,j_{l+1},j_{l+2}\}\big|\ge 2$ are true.

\item There are two different numbers $a$, $b$ in $\{0,1,\ldots, q-1\}$ satisfying one of
$$\begin{cases}
(M_{i_1,j_1},M_{i_2,j_2},\ldots,M_{i_k,j_k})=(a,b,a,b,\ldots,a,b)\\
(M_{i_1,j_1},M_{i_2,j_2},\ldots,M_{i_k,j_k})=(a,b,a,b,\ldots,a,b,a)
\end{cases}.$$
\end{enumerate} 
Denoting $M_{i_1,j_1}=M_{i_{k+1},j_{k+1}}$ and $M_{i_2,j_2}=M_{i_{k+2},j_{k+2}}$, if the sequence $$(M_{i_1,j_1},M_{i_2,j_2},\ldots,M_{i_{k+2},j_{k+2}})$$ also holds Conditions 1-3, then we call the sequence $S$ a \emph{cycle} or an $ab$-\emph{cycle}. 
\begin{remark}\label{wremark}
We may assume that a path $S$ in $M$ has at most two entries of $M$ in each row and column. Otherwise, there is a subsequence of $S$ that forms a cycle and does not contain $M_{i_1,j_1}$ and $M_{i_k,j_k}$. Removing the entries of this
subsequence from $S$ forms a path shorter than $S$.
\end{remark}  

By the definition of cycle, a matrix with a cycle cannot be a weak lonesum matrix. In addition, we can construct an infinite family of forbidden matrices when $q\ge 5$ by using cycles. Let $M^n$ be a $5$-ary $n\times n$-matrix defined by
$$\begin{cases}
M^n_{i,i}=0&\text{for}\ i\in\{1,2,\ldots,n\} \\
M^n_{i,i+1}=1&\text{for}\ i\in\{1,2,\ldots,n-1\}\\
M^n_{n,1}=1\\
M^n_{i,j}=2&\text{for}\ i\in\{2,3,\ldots,n-1\}\ \text{and}\ i\ge j+1\\
M^n_{i,j}=3&\text{for}\ i\in\{1,2,\ldots,n\}\ \text{and}\ i\leq j-2\\
M^n_{i,1}=4&\text{for}\ i\in\{2,3,\ldots,n-1\}
\end{cases}.$$
For example,
$$M^5 = \begin{pmatrix}
\mathbf{0}&\mathbf{1}&3&3&3\\
4&\mathbf{0}&\mathbf{1}&3&3\\
4&2&\mathbf{0}&\mathbf{1}&3\\
4&2&2&\mathbf{0}&\mathbf{1}\\
\mathbf{1}&2&2&2&\mathbf{0}\\
\end{pmatrix}.
$$
\begin{proposition}\label{iwlmatrix}
For $n\ge 3$ the matrix $M^n$ is not a weak lonesum matrix,
while each of its proper submatrices is a weakly lonesum matrix.
\end{proposition}
\begin{pf}
Interchanging $0$s and $1$s in $M^n$ does not change the row and column structures of $M^n$, thus $M^n$ is not a weak lonesum matrix.
However, the row and column structures of $M^n$
uniquely determine positions of $2$s, $3$s, and $4$s in $M^n$.
Therefore the remaining $2n$ entries can be filled by only two ways, 
and removing any row (resp. column) of $M^n$ yields a different column (resp. row) structure,
which means that every proper submatrix of $M^n$ is a weak lonesum matrix. 
\end{pf}
Proposition \ref{iwlmatrix} implies that we can generate an infinite sequence 
of $q$-ary forbidden matrices when $q\ge 5$. However, this is not true when $q\in\{3,4\}$. 
\begin{theorem}\label{iwlmatrix1}
Every quarternary matrix $M$ with a cycle of length at least $6$
contains a $2\times 2$, $2\times 3$, or $3\times 2$-matrices that are not weak lonesum matrices.
\end{theorem}
\begin{pf}
Suppose that a quarternary matrix $M$ with a minimal cycle of length $2n\ge 6$ does not satisfy the assumption of theorem.
By Remark \ref{wremark} we may assume that 
$M$ is an $n\times n$-matrix that contains a miniaml $01$-cycle of the form
$$(M_{0,0},M_{0,1},M_{1,1},M_{1,2},\ldots,M_{n-1,n},M_{n,n},M_{n,1}).$$
Since this is a minimal cycle, all other entires of $M$ should be either $2$ or $3$:
$$
\begin{pmatrix}
\mathbf{0}&\cdot&&\mathbf{1}&\\
&\cdot&\cdot&&\\
&&\cdot&\cdot&\\
&&&\mathbf{0}&\mathbf{1}\\
\mathbf{1}&&&&\mathbf{0}\\
\end{pmatrix},
\begin{pmatrix}
\cdot&\mathbf{1}&&\mathbf{0}&\\
&\mathbf{0}&\mathbf{1}&&\\
&&\mathbf{0}&\mathbf{1}&\\
&&&\cdot&\cdot\\
\cdot&&&&\cdot\\
\end{pmatrix} 
$$
Consider the top-left $3\times 3$ submatrix of $M$
$$
\begin{pmatrix}
\mathbf{0}&\mathbf{1}&c\\
a&\mathbf{0}&\mathbf{1}\\
&b&\mathbf{0}
\end{pmatrix}. 
$$
We know that $a$, $b$, and $c$ are either $2$ or $3$.
Moreover, if either $a=c$ or $b=c$, then $M$ contains a forbidden $2\times 3$
or $3\times 2$-matrices. Hence $a\neq c$ and $b\neq c$, which yields $a=b$.
Similarly, we conclude that
$$\begin{cases}
M_{2,1}=M_{3,2}=\ldots =M_{n,n-1}=M_{1,n}\\
M_{1,n}\neq M_{1,3}\\
M_{1,3}=M_{2,4}=\ldots =M_{n-2,n}=M_{n-1,1}=M_{n,2}\end{cases}.$$
Without loss of generality, we assume $M_{2,1}=2$.

Now we consider the element $M_{i,j}$ where $i\equiv j+2\ (\bmod\ n)$ and $j\in\{1,2,\ldots,n\}$. The matrix $M$ satisfies $M_{i,i}=M_{j,j}=0$ and
$M_{j,i}=3$, thus if $M_{i,j}=3$ then $M$ contains a forbidden
$2\times 2$-matrix. The only remaining case is $M_{i,j}=2$.
 
Similarly, if we consider $M_{j-3,j}$ for $j\in\{1,2,\ldots ,n\}$ on the basis of
$M_{j-3,j-2}=M_{j-1,j}=1$ and $M_{j-1,j-2}=2$, then we can conclude
that $M_{j-3,j}=3$ where we consider the indices of $M_{j-3,j}$ under modulo $n$.
Then the equalities $M_{j+3,j+3}=M_{j,j}=0$ and $M_{j,j+3}=3$ yield $M_{j+3,j}=2$. Likewise, $M_{j-4,j-3}=M_{j-1,j}=1$ and $M_{j-1,j-3}=2$ produce
that $M_{j-4,j}=3$. Continuing this argument provides
$$M_{j+k,j}=2\ \text{and}\ M_{j-k-1,j}=3$$ 
for $k\in\Big\{1,2,\ldots ,\Big\lfloor\frac{n-1}{2}\Big\rfloor\Big\}$. If we repeat this argument for 
$k=\Big\lfloor\frac{n-1}{2}\Big\rfloor+1$, then this yields a contradiction.
\end{pf}
\begin{remark}
Every binary non-lonesum matrix contains a cycle, thus Theorem 
\ref{iwlmatrix1} gives a third proof of
Theorem \ref{ryser}.
\end{remark}

In the rest of this subsection, we introduce some
examples of forbidden ternary matrices and say some informal words about their structures.

Two examples of forbidden ternary matrices are 
$$
T=\begin{pmatrix}
\mathbf{0}_1&\mathbf{1}_2&\mathbf{2}_0&0&0&0\\
1&1&\mathbf{0}_1&0&0&\mathbf{1}_0\\
1&1&\mathbf{1}_2&\mathbf{2}_0&\mathbf{0}_1&1\\
1&\mathbf{2}_1&2&2&\mathbf{1}_2&1\\
\mathbf{1}_2&2&2&2&\mathbf{2}_0&\mathbf{0}_1\\
\mathbf{2}_0&2&2&\mathbf{0}_2&0&0\\
\end{pmatrix}, T^{\prime}=
\begin{pmatrix}
1&\mathbf{0}_1&\mathbf{1}_2&\mathbf{2}_0&2&0&0&0&0\\
1&1&1&\mathbf{0}_2&\mathbf{2}_1&0&0&0&\mathbf{1}_0\\
1&1&1&2&\mathbf{1}_2&\mathbf{2}_0&\mathbf{0}_1&0&1\\
1&1&\mathbf{2}_1&2&2&2&\mathbf{1}_0&\mathbf{0}_2&1\\
\mathbf{1}_2&1&2&2&2&2&0&\mathbf{2}_0&\mathbf{0}_1\\
\mathbf{2}_1&\mathbf{1}_0&2&2&2&\mathbf{0}_2&0&0&0\\
\end{pmatrix}.
$$
The indices on entries in  $T$ (resp. $T^{\prime}$) indicate the entires of an alternative matrix with the same row and column structure vectors of $T$ (resp. $T^{\prime}$).
We can obtain more forbidden matrices by combining the following interchanges: 
$$\begin{cases}
\text{The $1$st column of $T$}\ \Longleftrightarrow\ 
\text{The $1$st and $2$nd columns of $T^{\prime}$}\\
\text{The $3$rd column of $T$}\ \Longleftrightarrow\
\text{The $4$th and $5$th columns of $T^{\prime}$}\\
\text{The $5$th column of $T$}\ \Longleftrightarrow\
\text{The $7$th and $8$th columns of $T^{\prime}$}
\end{cases}.$$

We can observe that some rows and columns in $T$ and $T^{\prime}$
contain triples of entries  with
values either $0_1$, $1_2$, $2_0$ or $0_2$, $1_0$, $2_1$
We call such triples $3$-\emph{trades}.
Every entry of a $3$-trade with the value $a_b$ is connected
by an $ab$-path with an entry of another $3$-trade, hence the value of the final entry is either $b_a$ (if the starting and
finishing $3$-trades are in the same row or column of a matrix)
or $a_b$ (otherwise). Every such a path in $T$ and $T^{\prime}$
has length $0$, $1$, or $2$. If we define the parity of a $3$-trade in such a way
that each path connects an even or odd numbers of $3$-trades,
then the incidence between paths and $3$-trades corresponds to
the incidence between edges and vertices of a bipartite cubic
graph. We can find a complete bipartite graph $K_{3,3}$ in $T$ or $T^{\prime}$ by this method. 

For a every
ternary weak non-lonesum matrix without cycles, 
we can find a similar structure corresponding to some other
bipartite cubic graph (in general, multiedges are allowed).
If there is an infinite sequence of ternary forbidden matrices,
then it corresponds to a sequence of cubic graphs. We can easily show that
the matrices with corresponding graphs having $2$ or
$4$ vertices contain a forbidden $2\times 2$, $2\times 3$, or
$3\times 2$-matrices. So, the matrices $T$ and $T^{\prime}$
are the minimal possible ones in some sense.
We also conjecture that the size $6\times 6$ is the smallest one for such matrices. 
A structure with six $3$-trades can be replaced by a $5\times 6$-matrix of the form
$$M=
\begin{pmatrix}
\mathbf{0}_1&&a&&&\mathbf{1}_0\\
\mathbf{1}_2&\mathbf{2}_0&\mathbf{0}_1&&&\\
\mathbf{2}_0&&b&\mathbf{0}_1&\mathbf{1}_2& \\
&\mathbf{0}_2&&\mathbf{1}_0&&\mathbf{2}_1\\
&&\mathbf{1}_0&&\mathbf{2}_1&\mathbf{0}_2\\
\end{pmatrix}.
$$ 
However, every way to fill the remaining entries of $M$ forms a forbidden
$2\times 2$ or $3\times 2$-submatrices
(for this, it is sufficient to consider the values of $a$ and $b$).

Using the concept of cycle, we have constructed an infinite 
sequence of $5$-ary forbidden matrices, 
which are also suitable for any $q$-ary case with $q\ge 5$.
When $q\leq 4$, we have shown
that a cycle always forms a trivial forbidden submatrix.
While nontrivial examples of ternary forbidden matrices exist,
to verify that infinite sequence of ternary or quarternary forbidden matrices exist remains open.

\section*{Acknowledgements}
We thank Prof. Kaneko for introduction of his results on poly-Bernoulli numbers and 
Sergey Avgustinovich for helpful discussions on weak lonesum matrices.

\bibliographystyle{elsarticle-num}

\bibliography{Poly-Bernoulli_numbers_and_lonesum_matrices}

\end{document}